\documentclass{amsart}

\usepackage[
    paperheight=11in,
    paperwidth=8.5in,
    outer=1.25in,
    inner=1.25in,
    bottom=1.4in,
    top=1.4in
]{geometry}
\usepackage{amsmath,amsxtra,amsthm,amsfonts,amssymb}
\usepackage{mathrsfs}
\usepackage[T1]{fontenc}

\DeclareSymbolFont{pploperators}{OT1}{ppl}{m}{n}
\SetSymbolFont{pploperators}{bold}{OT1}{ppl}{b}{n}
\makeatletter
\renewcommand{\operator@font}{\mathgroup\sympploperators}
\makeatother

\usepackage{graphicx,url,color,bbm,tikz-cd,enumitem}
\usepackage{microtype,xurl,needspace}
\usepackage[hidelinks,backref=page]{hyperref}
\hypersetup{
    colorlinks,
    citecolor=magenta,
    filecolor=magenta,
    linkcolor=blue,
    urlcolor=black,
    pdftitle={A tropical flag that cannot be completed}
}

\theoremstyle{definition}
\newtheorem{thm}{Theorem}[section]

\newtheorem{defn}[thm]{Definition}
\newtheorem{rem}[thm]{Remark}
\newtheorem{maintheorem}{Theorem}

\newcommand{\RR}{\mathbb R}

\DeclareMathOperator{\FlDr}{FlDr}
\DeclareMathOperator{\val}{val}
\numberwithin{equation}{section}

\title{Partial flags of tropical linear spaces may not complete}
\author{Christopher Eur}
\date{\today}

\begin{document}

\begin{abstract}
We construct two valuated matroids of ranks 2 and 4 on eight elements that form a valuated quotient, such that no rank 3 valuated matroid forms a valuated quotient with both simultaneously.
In other words, a partial flag of tropical linear spaces need not extend to a complete flag containing every intermediate dimension.
\end{abstract}

\maketitle

\vspace{-10pt}
\section{Introduction}

For a finite set $E$ and an integer $0\leq r \leq |E|$, let $\binom{E}{r}$ denote the set of $r$-subsets of $E$.
Write $\overline{\RR} := \RR \cup \{\infty\}$, with $\infty+a=\infty$.
We write $I\cup j$ for $I\cup\{j\}$ and $J\setminus j$ for $J\setminus\{j\}$.
Introduced by Dress and Wenzel \cite{DW} and interpreted as a \emph{tropical linear space} by Speyer \cite{Speyer}, a \emph{valuated matroid} of rank $r$ on $E$ is a function $\mu:\binom Er\to \overline{\RR}$, not identically $\infty$, such that for all subsets $I$ and $J$ of $E$ with $|I|=r-1$ and $|J|=r+1$, the minimum in
\[
 \min_{j\in J\setminus I}\bigl(\mu(I\cup j)+\mu(J\setminus j)\bigr)
\]
is attained at least twice (where $\infty$ is allowed to be the minimum).
We say that $\mu$ satisfies the \emph{tropical Pl\"ucker relations} in this case.
The subsets $B\in \binom{E}{r}$ such that $\mu(B) < \infty$ form the set of bases of a matroid, called the \emph{underlying matroid} of $\mu$.
A valuated matroid is said to be \emph{realizable} if, up to an additive constant, it is given by the valuations of the maximal minors of a matrix over a valued field.
For background on matroids and valuated matroids, we refer to \cite[Chapter~4]{MS}.

\medskip
Inclusions of tropical linear spaces are described in terms of their valuated matroids as follows.

\begin{defn}\label{def:flag}
Let $f$ and $g$ be valuated matroids on $E$ of ranks $r\leq s$.
We say that $f$ is a \emph{valuated quotient} of $g$, denoted $f\twoheadleftarrow g$, if for all subsets $I$ and $J$ of $E$ with $|I|=r-1$ and $|J|=s+1$, the minimum in
\begin{equation}\label{eq:incidence}
 \min_{j\in J\setminus I}\bigl(f(I\cup j)+g(J\setminus j)\bigr)
\end{equation}
is attained at least twice (again allowing $\infty$ as the minimum).
We say that $(f,g)$ satisfies the \emph{tropical incidence relations} in this case.
A \emph{valuated flag matroid} is a sequence $(f_1,\ldots,f_k)$ of valuated matroids with $f_i\twoheadleftarrow f_j$ for every $i<j$. The $f_i$ are called its \emph{constituents}.
\end{defn}

Brandt--Eur--Zhang \cite[Section~3.2 and Theorem~4.3.1]{BEZ} showed that a rank $r$ valuated matroid $f$ determines a (projective) tropical linear space $L(f)$ of dimension $r-1$, and that $f\twoheadleftarrow g$ if and only if $L(f)\subseteq L(g)$.
In this short note, we record an example showing that, unlike flags of linear subspaces over a field, a partial flag of tropical linear spaces may not admit a completion to a complete flag.

\begin{maintheorem}\label{thm:main}
There are valuated matroids $p$ and $q$ on eight elements of ranks $2$ and $4$ (respectively), each realizable and with uniform underlying matroids, such that $p\twoheadleftarrow q$ but no rank $3$ valuated matroid $h$ satisfies $p\twoheadleftarrow h\twoheadleftarrow q$.
\end{maintheorem}

In the notation of \cite{BEZ}, the forgetful map
\[
 \FlDr(2,3,4;8)\longrightarrow\FlDr(2,4;8)
\]
is therefore not surjective.
This answers negatively Wang's completion question \cite[\textup{(P3)} and Section~6.2]{Wang}.
The problem of inserting constituents of missing ranks also appears in \cite[Section~5.1]{BLS}.
A stronger question asks whether the two valuated matroids in a valuated quotient can be obtained by deleting and contracting the same set in a larger valuated matroid \cite[Remark~2.24]{JL}, \cite[Remark~5.8]{BLS}.
Since such a larger valuated matroid induces a completion, the example here also rules out such an extension (Remark~\ref{rem:extension}).

\section{The eight-element example}
We construct $p$ and $q$ in Theorem~\ref{thm:main} as follows.
Partition $E=[8] = \{1, \dotsc, 8\}$ into four pairs
\[
 A=\{1,2\},\qquad B=\{3,4\},\qquad C=\{5,6\},\qquad D=\{7,8\},
\]
and call $A\cup B$ and $C\cup D$ the two \emph{halves}.
For $S\subseteq E$, let $d(S)$ be the number of these pairs wholly contained in $S$.
We omit braces and commas when listing the elements of a set, writing $ij$ for $\{i,j\}$, and so on.
Define
\[
 p(ij)=\begin{cases}
 2,&i,j\text{ belong to the same pair},\\
 1,&i,j\text{ belong to different pairs in the same half},\\
 0,&i,j\text{ belong to opposite halves},
 \end{cases}
 \qquad\text{and}\quad q(S)=d(S)\quad\left(S\in\binom E4\right).
\]
Thus $p$ records two levels of clustering, first into halves and then into pairs, whereas $q$ records only the four pairs.
Let $K$ be any valued field such that 2 and 3 are units in $K$, and let $t\in K$ be an element with valuation 1.

\begin{proof}[Proof of Theorem~\ref{thm:main}]
Label the elements lexicographically by $(a,b,c)\in\{0,1\}^3$, so that $a$ specifies a half and $(a,b)$ specifies a pair.
Define
\[
 z_{abc}=a+bt+ct^2\qquad\text{and}\qquad w_{abc}=2a+b+ct.
\]
The matrix with columns $(1,z_i)^{\mathsf T}$ realizes $p$.
The matrix with columns $(1,w_i,w_i^2,w_i^3)^{\mathsf T}$ realizes $q$: a maximal minor is a Vandermonde product, and $\val(w_j-w_i)$ is $1$ within a pair and $0$ otherwise.

To check $p\twoheadleftarrow q$, fix $i\in E$ and $J\in\binom{E}{5}$.
For $j\in E$, we write $j'$ for the other element of the pair containing $j$.
Since $d(J\setminus j)=d(J)-\mathbf1_{\{j'\in J\}}$, subtracting $d(J)$ from each term of \eqref{eq:incidence}, we obtain for each $j\in J\setminus i$ the entry
\[
 p(ij)-\mathbf1_{\{j'\in J\}}.
\]
We verify that the minimum occurs at least twice among these entries for all choices of $i$ and $J$.
If $J$ contains a whole pair in the half opposite $i$, its two elements give the minimum $-1$.
Otherwise $J$ has one or two elements in that half.
If there are two, which do not form a pair, their entries are both 0 and all other entries are nonnegative.
If there is one, $J$ contains all four elements of $i$'s half; the pair not containing $i$ then gives two entries equal to 0, again minimal.
This proves $p\twoheadleftarrow q$.

We now exclude a rank $3$ intermediate valuated matroid $h$ satisfying $p\twoheadleftarrow h\twoheadleftarrow q$.
Let us divide the triples in $E$ into three types:
\[
\begin{array}{c|l}
X&\text{a whole pair and a singleton in the same half},\\
Y&\text{a whole pair and a singleton in the opposite half},\\
Z&\text{one element from each of three distinct pairs}.
\end{array}
\]
Let $x,y,z\in\overline{\RR}$ be the respective minima of $h$ on these types.
We will prove
\begin{equation}
 x\geq\min\{y+1,z+2\},\qquad z\geq x-1,\qquad y\geq z+1.
\end{equation}
This gives the desired nonexistence of the $h$.
If $x=\infty$, the last two inequalities force $z=y=\infty$, contradicting the requirement that $h$ have a finite coordinate.
If $x<\infty$, they give $y+1\geq z+2\geq x+1$, and the first inequality forces $x\geq x+1$, again a contradiction.

For the proof of the above inequalities, for each minimum $x,y,z$, we choose a triple attaining it in the casework below.
We continue to use $j'$ for the element forming a pair with $j\in E$.
\begin{enumerate}
\item[($X$)] For a type-$X$ triple $aa'b$, choose $c$ in the opposite half.
The tropical incidence relation for $p\twoheadleftarrow h$ with $I=a$ and $J=aa'bc$ is
\[
 \min\{h(aa'b),\ 1+h(aa'c),\ 2+h(abc)\}.
\]
Since the first term $h(aa'b)=x$ cannot be uniquely minimal, and since the other two terms are at least $1+y$ and $2+z$, we find $x\geq\min\{y+1,z+2\}$.
\item[($Y$)] For a type-$Y$ triple $aa'c$, choose $b$ in the other pair of $a$'s half and $d$ in the other pair of $c$'s half.
The tropical incidence relation for $h\twoheadleftarrow q$ with $I=ac$ and $J=aa'bcd$ is
\[
 \min\{h(aa'c),\ 1+h(abc),\ 1+h(acd)\}.
\]
Since $h(aa'c)=y$, the other two terms, both from type $Z$, give $y\geq z+1$.
\item[($Z$)] For a type-$Z$ triple $abc$, we may assume that $a,b$ belong to the same half.
The tropical incidence relation for $h\twoheadleftarrow q$ with $I=ab$ and $J=aa'bb'c$ is
\[
 \min\{1+h(aa'b),\ 1+h(abb'),\ 2+h(abc)\}.
\]
Since $h(abc)=z$, the first two terms, both from type $X$, give $z+2\geq x+1$.\qedhere
\end{enumerate}
\end{proof}

We note that the proof uses only the tropical incidence relations to rule out the existence of $h$; that is, we did not need to consider the tropical Pl\"ucker relations for $h$.

\begin{rem}[Common extensions]\label{rem:extension}
There is no valuated matroid $R$ on $E\sqcup T$, for any finite set $T$ disjoint from $E$, such that $R\setminus T=q$ and $R/T=p$.
Indeed, replacing deletions by contractions one at a time gives a flag of valuated quotients (see \cite[Corollary~4.3.2(1), (3)]{BEZ} or \cite[Theorem~2.23]{JL}).
Its ranks decrease from $4$ to $2$ by steps of at most one, supplying the rank $3$ intermediate excluded by Theorem~\ref{thm:main}.
\end{rem}

\subsection*{AI disclosure}
This example was found entirely by ChatGPT-6 Astra on Work.  The author verified the argument and polished the exposition.

\small
\bibliographystyle{alpha}
\bibliography{Eur_NoncompletableFlag}

\begin{thebibliography}{Wan24}

\bibitem[BEZ21]{BEZ}
Madeline Brandt, Christopher Eur, and Leon Zhang.
\newblock Tropical flag varieties.
\newblock {\em Advances in Mathematics}, 384:107695, 2021.
\newblock \url{https://doi.org/10.1016/j.aim.2021.107695}.

\bibitem[BLS26]{BLS}
Marie-Charlotte Brandenburg, Georg Loho, and Ben Smith.
\newblock Quotients of {$M$}-convex sets and {$M$}-convex functions.
\newblock {\em Combinatorial Theory}, 6(1):Paper No.~14, 53 pp., 2026.
\newblock \url{https://doi.org/10.5070/C66165701}.

\bibitem[DW92]{DW}
Andreas W.~M. Dress and Walter Wenzel.
\newblock Valuated matroids.
\newblock {\em Advances in Mathematics}, 93(2):214--250, 1992.
\newblock \url{https://doi.org/10.1016/0001-8708(92)90028-J}.

\bibitem[JL24]{JL}
Manoel Jarra and Oliver Lorscheid.
\newblock Flag matroids with coefficients.
\newblock {\em Advances in Mathematics}, 436:109396, 2024.
\newblock \url{https://doi.org/10.1016/j.aim.2023.109396}.

\bibitem[MS15]{MS}
Diane Maclagan and Bernd Sturmfels.
\newblock {\em Introduction to Tropical Geometry}, volume 161 of {\em Graduate
  Studies in Mathematics}.
\newblock American Mathematical Society, Providence, RI, 2015.
\newblock \url{https://doi.org/10.1090/gsm/161}.

\bibitem[Spe08]{Speyer}
David~E. Speyer.
\newblock Tropical linear spaces.
\newblock {\em SIAM Journal on Discrete Mathematics}, 22(4):1527--1558, 2008.
\newblock \url{https://doi.org/10.1137/080716219}.

\bibitem[Wan24]{Wang}
Jidong Wang.
\newblock Lorentzian polynomials and the incidence geometry of tropical linear
  spaces.
\newblock arXiv:2412.12059, 2024.
\newblock Version~2, 2026. \url{https://arxiv.org/abs/2412.12059v2}.

\end{thebibliography}

\end{document}